\documentclass[journal]{IEEEtran}

\usepackage{amsfonts,amssymb,amsbsy,textcomp,marvosym,picins,amsmath,caption,threeparttable,amsthm,subfigure}
\usepackage{eurosym,mathrsfs,fancyhdr,CJK,multicol,graphics,indentfirst,color,bm,upgreek,booktabs,graphicx}
\usepackage{amsmath,amstext,amssymb}
\usepackage{latexsym}
\usepackage{authblk}
\usepackage{cite}
\usepackage{acronym}
\usepackage{color}
\usepackage{epsfig}
\usepackage{epstopdf}
\usepackage{placeins}
\usepackage{subfig}
\usepackage{lipsum,widetext}

\usepackage{algorithm}
\usepackage{algpseudocode}
\usepackage{pifont}

\usepackage{csquotes}
\usepackage{wrapfig}
\begin{document}
\pagenumbering{gobble}

\title{Dynamic Constraint-based Influence Framework and its Application in Stochastic Modeling of Load Balancing}

\author[1]{Ehsan Siavashi}
\author[2]{Mahshid Rahnamay-Naeini}
\affil[1]{Department of Computer Science, Texas Tech University, Lubbock, TX, USA}
\affil[2]{Department of Electrical Engineering, University of South Florida, FL, USA}
\renewcommand\Authands{ and }

% The paper headers
%\markboth{IEEE Transactions on Power Systems}%
%{Shell \MakeLowercase{\textit{et al.}}: Bare Demo of IEEEtran.cls for Journals}

% make the title area
\maketitle

\begin{abstract}  Components connected over a network influence each other and interact in various ways. Examples of such systems are networks of computing nodes, which the nodes interact by exchanging workload, for instance, for load balancing purposes. In this paper, we first study the Influence Model, a networked Markov chain framework, for modeling network interactions and discuss two key limitations of this model, which cause it to fall short in modeling constraint-based and dynamic interactions in networks. Next, we propose the Dynamic and Constraint-based Influence Model (DCIM) to alleviate the limitations. The DCIM extends the application of the Influence Model to more general network interaction scenarios. In this paper, the proposed DCIM is successfully applied to stochastic modeling of load balancing in networks of computing nodes allowing for prediction of the load distribution in the system, which is a novel application for the Influence Model. The DCIM is further used to identify the optimum workload distribution policy for load balancing in networked computing systems.
\end{abstract}


\IEEEpeerreviewmaketitle

\section{Introduction}
Various networks in real-world, such as communication networks, power grids and social networks consist of connected components that interact with each other due to various reasons such as exchange of information or workload, physics- and engineered-based interactions and social rules. In the past decade, a large body of work has been emerged in modeling and analyzing the interaction of components over networks with applications in various domains including studying epidemics, information diffusion and spreading phenomenon (e.g., disease epidemics and cascading failures) [1-4]. Some of such interaction models are customized for specific problems and specific systems, such as cascading failures in power grids [4]. A common property of such network interaction models is the focus on simple component interactions to ensure the scalability and tractability of the models.
However, there are real-world scenarios of network interactions, such as interactions among the components in a network of computing elements, that do not quite fit into any of the available models due to constraints and rules affecting the interactions. The dynamic and constraint-based network interactions are natural attributes of many systems as in real-world the network components do not interact in a static and uniform manner but instead the magnitude and form of interactions vary in time and due to the statues of the components and the rules governing the system.

%Furthermore, the internal dynamics of a node could change due to the state of the node as well as its interactions with other nodes. Note that in this paper we use the terms node and network component interchangeably.

In order to enable dynamic constraint/rule-based network interactions, we present a model named {\it Dynamic and Constraint-based Influence Model (DCIM)} in which: (1) the internal dynamics of components can vary in time, for instance, based on their past interactions, and (2) depending on the state of the components the influences among components may get activated or deactivated (rule-based interactions).
The DCIM is an extension of the mathematical Influence Model (IM), a networked Markov chain framework, first presented in [5] and [6]. Specifically, IM associates a Markov chain (MC) with each node (note that in this paper we use the terms node and network component interchangeably) while the nodes interact with each other over an underlying network topology by affecting the transition probabilities of their individual MCs. We will argue that the original IM has limitations for modeling nodes with {\it varying internal dynamics} (note that we emphasis on varying internal dynamics not varying internal states) and constraint/rule-based interactions. The proposed DCIM alleviates these limitations and extends the application of the IM to more general interaction scenarios as well as spreading phenomena in networked systems. Although not all the nice mathematically tractable results of the original IM will hold for the extended DCIM as discussed in the paper, we will show that the extended IM enables various modeling and numerical analysis for more complex scenarios.
As an example of such scenarios, we present a novel application of the influence-based network interactions by using the DCIM to develop a stochastic load-balancing and distribution prediction model for a network of computing nodes. The DCIM also enables us to identify the best workload distribution policy or rule to achieve a desirable load-balancing distribution based on the prediction model. Specifically, the DCIM is used for identifying the optimum workload distribution policy to achieve a balanced computing network. The DCIM is adequately general to be applied to other network interaction scenarios such as disease epidemic studies, particularly for diseases that spread based on dynamic and constraint-based interactions, in order to provide more accurate epidemic models. More details and analysis of the proposed DCIM can be found in [7].

The rest of this paper is organized as following. We present an overview of the general IM in Section~II and discuss the novel application of stochastic load balancing for network of computing nodes in the language of IM as an example. Section~III is devoted to discussion about the limitations of the general IM. In Section~IV, we present the DCIM as a solution to the limitations introduced in Section~III and discuss the steady state behavior of the DCIM. In Section~V, we present the application of the DCIM to stochastic load balancing for network of computing nodes and present prediction results on the load distribution based on various load balancing policies and rules and different underlying topology of the system. We also present an algorithm to identify the optimum load balancing policy for the network of computing nodes. A review of related works is presented in Section~VI. Finally, in Section~VII, we present a summary of our conclusions.

%In this thesis, we will introduce an extension of IM, which allows modeling much more general interactions in networks. We have applied our proposed extended IM to a load balancing problem in a network of computing nodes. The model not only allows modeling the network interactions, but also allows identifying the best load balancing policy given the current state of the network. We have also applied and compared the application of the proposed model to a disease propagation model with and without immunization.
%
%\begin{figure*}
%  \centering \includegraphics[width=3 in]{DCNetwork}
%  \caption{Interdependent Markov chain model for coupled stochastic dynamics of operator state and infrastructure behavior.}
%\end{figure*}

\section{A Review of Influence Model (IM)}\label{SEC:Review}
\subsection{The general IM}
In this section, we review the germane aspects of IM as introduced in [5, 6] and also introduce the terminology and notations that will be used in this paper. As briefly mentioned in the Introduction Section, IM is a networked MC framework, which associates a MC with each network node and the nodes interact with each other based on an underlying network topology. In this model, the state evolution of each node depends on its internal MC as well as the state of its neighbors and their influences on the node (see Fig. 1 for an example of the influence network).

The IM considers a weighted and directed graph of $n$ nodes. The network topology and the weights of the links specify which nodes influence each other and what the strength of the influences are. The concept of these influences will be clarified with an example in the next subsection. The internal stochastic dynamics of each node are captured by a MC. The state space of the internal MC specifies the possible states that a node may attain. In this paper, we assume that the nodes of the network are homogeneous and have the same internal MC; however, the IM and the proposed DCIM in this paper are not limited to the homogeneous nodes and the results can be extended to heterogenous network nodes.
%The nodes may be heterogeneous and there is no restriction on the networks topology. Relax the homogenous assumption
The status of node $i$ at time $t$ is denoted by $s_i[t]$, a vector of length $m$, where $m$ is the number of possible states for the node (i.e., $m$ represents the cardinality of the state space of the MC representing the node's dynamics). At each time, all the elements of $s_i[t]$ are 0 except for the one which corresponds to the current state of the node (with value 1). We denote the $j$-th entry in the $s_i[t]$ by $s_{ij}[t]$. The status of the whole network encompassing all the nodes can be represented by concatenating individual state vectors into a vector $\mathbf{S}[t]=(s_1[t] s_2[t] ... s_n[t])$ of length $m\times n$. The transition probability matrix of the internal MC for node $i$ is $\mathbf{A}_{ii}$, which is an $m\times m$ row stochastic matrix.

In addition to the individual internal MCs for the nodes of the network, the influences between the nodes are also critical in characterizing and analyzing the behavior of the whole system. The influences enable capturing impacts of the state of a node on the state evolution of its neighbors. The influences of the network is captured by the influence matrix $\mathbf{D}$, which the entry $d_{ij}$ is a number between 0 and 1 representing the amount of influence that node $i$ receives from node $j$. The larger the $d_{ij}$ is the more influence the node $i$ receives from node $j$; with the two extreme cases being $d_{ij}=0$ meaning that node $i$ does not receive any influence from node $j$ and $d_{ij}=1$ meaning that the next state of node $i$ deterministically depends on the state of node $j$.  Note that receiving influence from a node itself, i.e., $d_{ii}$, specifies how much the state evolution of a node depends on its internal MC.  The total influence that a node receives should add up to unity i.e., $\sum^{n}_{j=1}d_{ij}=1$, and therefore, matrix $\mathbf{D}$ is a row stochastic matrix too.

The influence matrix $\mathbf{D}$ merely specifies how much two nodes influence each other. In order to specify how the states of the nodes will change due to the influences, we also need state-transition matrices $\mathbf{A}_{ij}$, which capture the probabilities of transiting to various states due to the state of the influencing node. Matrix $\mathbf{A}_{ii}$ represents the special case of self-influence, which is described by the internal MC of the node. Note that the $\mathbf{A}_{ij}$ matrices are row stochastic. In the general IM [1], the collective influences among the nodes in the network is summarized in the total influence matrix $\mathbf{H}$ define as:
\begin{equation}
  \mathbf{H}=\mathbf{D}^{\prime}\otimes\{\mathbf{A}_{ij}\}=
 \begin{pmatrix}
  d^{\prime}_{11}A_{11} & \cdots & d^{\prime}_{1n}A_{1n}\\
  \vdots  & \ddots & \vdots  \\
  d^{\prime}_{n1}A_{n1} & \cdots & d^{\prime}_{nn}A_{nn}
 \end{pmatrix},
\end{equation}
where $\mathbf{D}^{\prime}$ is the transpose of the matrix $\mathbf{D}$ and $\otimes$ is the generalized Kronecker multiplication of matrices [5]. Finally, based on the the total influence matrix $\mathbf{H}$ the evolution equation of the model is defined as
\begin{equation}
  \mathbf{p}[t+1]=\mathbf{S}[t]\mathbf{H},	
\end{equation}
where vector $\mathbf{p}[t+1]$ describes the probability of various states for all the nodes in the network in the next time step.

Steady state analysis of IM has some similarities with that of MCs.
%It follows from the definitions that the probability of node $i$ being in state $s_i$ at time $t$ is equal to $E(s_i[t])$, where $E(s_i[t])$ is the expected value of $s_i[t]$.
Specifically, given the initial state of the network, $\mathbf{S}[0]$, we have
\begin{equation}
  E(\mathbf{S}[t+1]|\mathbf{S}[0])=\mathbf{S}[0]\mathbf{H}^{t+1}, 	
\end{equation}
where $E(.)$ represents the expected value.
Therefore, having the initial state of the network and matrix $\mathbf{H}$ is sufficient for calculating the status probability vector, which can help determining the asymptotic probability of occupancy of the states for the system. According to [5], the steady state of the system in this case may or may not be independent of the initial state. For analyzing the steady-state of the network, we need to know if $\lim_{t \to \infty} \mathbf{S}[0]\mathbf{H}^t$ converges to any specific probability vector. It is shown in [5] that $\lim_{t \to \infty} \mathbf{S}[0]\mathbf{H}^t$ exists if $\mathbf{H}$ has an eigenvalue at 1 that dominates all of its other eigenvalues. For further reading on IM and its properties refer to [5] and [6].
In order to clarify the concepts of IM, we introduce an example application in the next subsection.
%
%The influence framework allows modeling of spreading phenomena in heterogeneous networks. In particular, the heterogeneity of network nodes can be due to different types of nodes in the network, which can attain various states and respond differently to influences and interactions with other components.
\begin{figure*}[tb!]
 	\centering
 	\includegraphics[width=4.0 in]{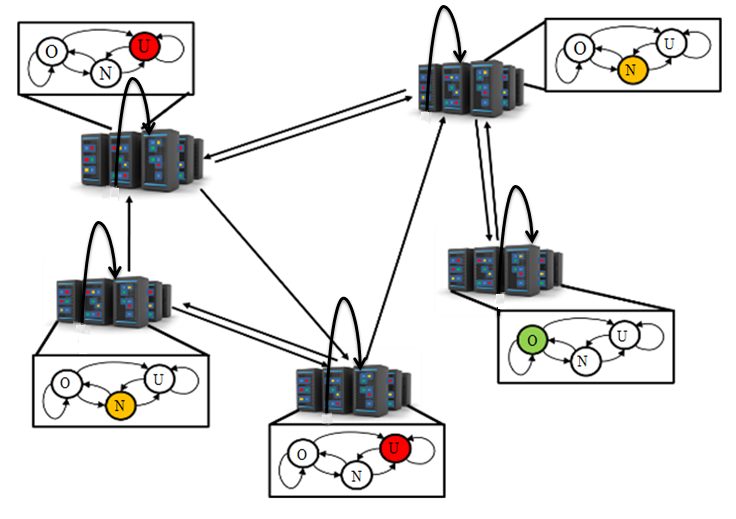}
 	\caption{A network of computing nodes modeled by the IM with internal MCs of the nodes and influences based on the topology of the network.}
 	\label{FIG:SYSM}
 \end{figure*}

% \begin{figure*}[!htb]
%\centering
%  \subfigure[]{
%    \includegraphics[width=6cm,height=2cm]{photo.eps}}
%  \subfigure[]{
%    \includegraphics[width=6cm,height=2cm]{photo.eps}}
%  \caption{Example for inserting a two-column wide figure. (a) Title of sub-figure (a). (b) Title of sub-figure (b).}
%\end{figure*}

\subsection{Modeling a Network of Computing Nodes with IM}
In this section, we introduce the application of the IM in modeling load distribution in a network of computing nodes.
Recently, a large body of work has been emerged to design and investigate the load-balancing mechanisms (with increasing performance and reducing energy consumption goals) in networks of computing elements (e.g., networks of data centers) [8, 9, 10]. Similar problems have been also studied in load balancing for multi-domain communication networks [11, 12] for network traffic instead of workloads/tasks in computing networks. The idea behind load-balancing mechanisms in general is to transfer workload among nodes to better utilize the computing resources by, for example, sending the workload from nodes, which are overloaded, to the nodes, which are underloaded, to improve the performance and the response time, or sending workloads from nodes with expensive energy resources to cheaper and greener resources. The exchange of workloads is easy to achieve with the help of advances in virtual machines (VMs) and VM migration techniques. The network of computing nodes in general is a highly stochastic environment with stochastic workload arrival, processing and exchange of workloads. In this subsection, we use the IM to model the stochastic interactions among the nodes in such networks.
This model will help clarifying the concepts of the IM in a physical setting. However, as will be discussed in the next section, the IM is not capable of modeling this system accurately by considering the load balancing rules. In Section VI, we will revisit this model and show that our proposed extension of the original IM will enables us to provide a more accurate model for this system by capturing the load balancing rules and dynamic constraint-based interactions.

An example of a network of computing nodes is depicted in Fig.~1 in which the distributed computing elements interact with each other over a network topology. Note that in this application the topology of the underlying network could be the physical topology of the communication network or a logical topology overlaid on top of the physical topology for distributing workload among nodes. In a coarsely quantized setting and according to the size of the workload, we can describe the state of the nodes by, for example, overload ($O$), normal ($N$) or underload ($U$), i.e., the state space of the internal MC for each computing element is the set $\{O, N, U\}$. The stochastic transitions among the internal states of nodes depend on the stochastic workload processing and the external workload arrivals to the node. In other words, the state of a node may change not only because of its internal stochastic dynamics (due to processing the workload) but also due to receiving influences (workloads) from other nodes.

\begin{figure*}[tbh!]
 	\centering
 	\includegraphics[width=4 in]{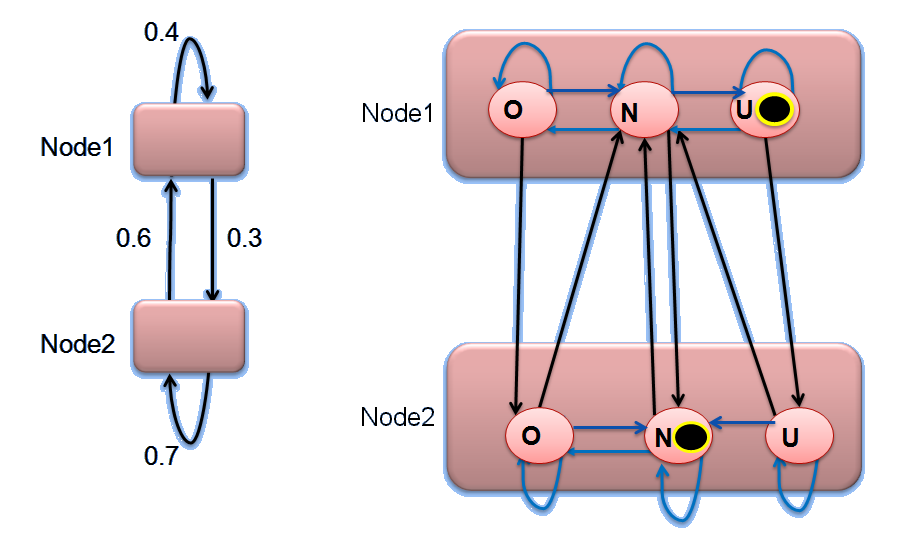}
 	\caption{An example of a network of computing nodes with two nodes. (a) The elements of the influence matrix for interactions between the nodes, i.e., matrix $\mathbf{D}$, and (b) the elements of the state-transition matrices and the internal MCs capturing the dynamics of nodes, i.e., matrices $\mathbf{A}_{ij}$.}
 	\label{FIG:SYSM2}
 \end{figure*}

Let us explain the concepts of IM by focusing on a simple scenario in which the network is consisted of two computing nodes as depicted in Fig. 2. We will expand this example to more than two nodes with arbitrary topology in Section VI. Figure~2-a shows the elements of the influence matrix, $\mathbf{D}$, for this network. For instance, node 1 receives $60\%$ of its total influences from node 2 (by receiving workload) and the rest from its internal MC. Figure~2-b depicts the details of dynamics inside each node by showing the elements of state-transition matrices, i.e., $\mathbf{A}_{ij}$s, including individual MCs. In particular, there are two types of arrows in Fig.~2-b: (1) arrows that are inside the square and represent the internal MC transitions of the nodes, i.e., matrices $\mathbf{A}_{ii}$, and (2) arrows that cross the boundaries of the sites represent the state-transitions, i.e., matrices $\mathbf{A}_{ij}$. The black dot inside the nodes indicates the current state of the node. For example, the current states of nodes 1 and 2 in Fig.~2-b are $s_1[t]=(0, 0, 1)$ and $s_2[t]=(0, 1, 0)$, respectively.
 According to the internal MCs, the probability of occupying a state by a dot is positive if and only if it is currently one step (an arrow with a positive weight) away from any dot.
 For instance, solely based on the MC of node 2, the probability transiting to state $U$ is 0; however, because of the influences in the IM node 2 can go to the state $U$ in the next step although there is not direct transition link inside the MC. In a network of computing nodes if a node is currently in an underload state and there is no chance for it to go to the overload state in the next step based on its internal MC; however, it still may go to the overload state in the next time step if other nodes sent a huge workload on the node, say due to failure of a node.
Note that the workload distribution and transfer among the nodes depend on the load balancing rules and the state of the nodes in the network. We elaborate on load distribution rules in Section VI.
\vspace{-0.2cm}
\section{Limitations of the Original Influence Model}
In this work, we have identified two key limitations of the original IM, which are discussed next.
\subsubsection{IM Lacks the Capability of Capturing Constraint-Based Influences}
The IM in its current form cannot model scenarios in which the influence of node $i$ on node $j$ can get activated or deactivated depending on certain constraints or rules. Let us explain this scenario with an example. Recall that, in the example of network of computing elements, influencing a node means sending workload to the node. Suppose that the load balancing policy in the network suggests that if a computing node is overloaded (in state $O$) and its (physical/logical) neighbor according to the network topology is underloaded (in state $U$) then the node should influence (send workload) to the neighbor. This means that a node will not influence other nodes if it is in normal (in $N$ state) or underloaded (in $U$ state). Capturing the latter cases in the state-transition matrices using IM will lead to all zero rows, which violates the row stochastic property of the $\mathbf{A}_{ij}$ matrices. Therefore, the model cannot capture such scenarios with rule-based interactions depending on the state of the nodes.

Moreover, in the mentioned load-balancing rules, not only the state of the influencing node but also the state of the node, which is receiving influence, is important. For example, a node should only influence a neighbor if the neighbor is underloaded (in $U$ state). Similarly, as the state-transition matrices for capturing the transitions due to influences are static in the original IM, one cannot activate or deactivate the influences depending on the state of the nodes to capture these scenarios.

\subsubsection{IM Lacks the Capability of Capturing Changes in the Internal Dynamics of MCs}
In many real world networks, the dynamics of state transitions of a node (modeled by $\mathbf{A}_{ii}$ in IM) may change in time or due to various reasons. For instance, there are real-world scenarios that the action, in which  node $i$ influences node $j$, not only has an effect on node $j$, but also has a side effect on node $i$ itself. This happens, for example, when a computing node sends workload to another node. The action of transferring workload not only increases the likelihood of transiting to a state with larger workload for the receiving node, but also increases the probability of the influencer node to transit to a less loaded state in the next time step. As the state transition matrices in the general IM are static and do not change depending on the influences among nodes, these scenarios cannot be captured in the model. In DCIM, we allow dynamic change of the transition matrices for internal MCs.
\vspace{-0.2cm}
\section{Dynamic and Constraint-based Influence Model (DCIM)}
In this section, we introduce Dynamic and Constraint-based Influence Model (DCIM). The DCIM addresses the first problem discussed in Section~V by defining a constraint matrix $\mathbf{C}$, where the entry $c_{ij}$ for $i,j \leq n$ is a binary variable specifying whether node $i$ gets influenced by node $j$ or not. Specifically, $c_{ij}= 1$ indicates that node $i$ gets influences by node $j$ and $c_{ij}= 0$ indicates otherwise. We assume that each node always influences itself (i.e., $c_{ii} = 1$ for all $0 < i \leq n$). One can define the value of $c_{ij}$ according to boolean logic to capture the rules and policies for interactions in the network. In other words, $c_{ij}$s are functions of the state of the nodes. For instance, in a network of computing elements if the load balancing policy dictates that node 1 gets influenced by node 2 only if node 2 is overloaded and node 1 is underloaded, then we can define $c_{12}=s_{13}s_{21}$, where $s_{ij}$ is the $j$-th entry of the state vector of the node $i$ as explained in Section~II.

Now, using the constraint matrix $\mathbf{C}$ and the influence matrix $\mathbf{D}$ as introduced in Section~II, we define the constraint-based influence matrix for DCIM, denoted by $\mathbf{E}$, as
\begin{equation}
  \mathbf{E}=\mathbf{D} \circ \mathbf{C}+\mathbf{I} \circ(\mathbf{D}\times (\mathbf{1}-\mathbf{C}^{\prime})),
  \label{EQ:E}
\end{equation}
where $\circ$ is the Hadamard product (aka entrywise product), $\mathbf{1}$ is the matrix with all elements equal to 1 and $\mathbf{C}^{\prime}$ is the transpose of matrix $\mathbf{C}$. Based on the above definitions and assumptions, matrix $\mathbf{E}$ will be of the following form:

\vspace{-0.3cm}
\begin{equation}
\mathbf{E}\!\!=\!\!
 \begin{pmatrix}
  d_{11}+\sum_{i\neq 1}^n d_{1i}(1-c_{1i}) & d_{12}c_{12} & \!\!\!\!\!\!\cdots & d_{1n}c_{1n} \\
  d_{21}c_{21} & \!\!\!\!\!\!\!\!\!\!\!\!\!\!\!\!\!\!d_{22}+\sum_{i\neq 2}^n d_{2i}(1-c_{2i}) & \!\!\!\!\!\!\cdots & d_{2nc_{2n}} \\
  \vdots  & \vdots  & \!\!\!\!\!\!\ddots & \vdots  \\
  d_{n1}c_{n1} & d_{n2}c_{n2} & \!\!\!\!\!\!\cdots &
 \end{pmatrix}
\end{equation}
Note that in this model, matrix $\mathbf{E}$ plays a similar role to that of matrix $\mathbf{D}$ in the original IM. More specifically, we have
 \begin{equation}
  \mathbf{H}=\mathbf{E}^{\prime}\otimes\{\mathbf{A}_{ij}\},
  \label{EQ:NewH}
\end{equation}
and $\mathbf{p}[t+1]=\mathbf{S}[t]\mathbf{H}$.

Next, we investigate some of the properties of this model.

\emph{Proposition I:} The total constraint-based influence matrix $\mathbf{E}$ is row stochastic.
\begin{proof}
By expanding the summation of a row, it is easy to observe that regardless of the values of $c_{ij}s$ we have $\sum_j^ne_{ij}=\sum_j^nd_{ij}=1$.
\end{proof}

Now, let $\mathbf{S}_i[t]$ and $\mathbf{p}_i[t]$ be the sub-vectors of $\mathbf{S}[t]$ and $\mathbf{p}[t]$, which, respectively, represent the state and probability vectors for just node $i$.
For each node $i$, $(0<i\leq n)$ we have
\begin{equation}
\begin{aligned}
\mathbf{p}_i[t+1] &= (d_{ii}+\sum_{j \neq i}^n d_{ij}(1-c_{ij}))\mathbf{S}_i[t]\mathbf{A}_{ii}\\
&+ (\sum_{j\neq i}^n d_{ij}c_{ij}S_j[t])\mathbf{A}_{ji}.
\end{aligned}
\label{EQ:expand}
\end{equation}
Based on (\ref{EQ:expand}) it is easy to derive the following propositions.
\emph{Proposition II:} For any node $i$ that receives no influence from other nodes in the network, node $i$'s next state depends only on its internal MC and its current state.
\begin{proof}
By setting $c_{ij}=0$ for $j\neq i$ in (\ref{EQ:expand}), we get $\mathbf{p}_i[t+1]=\mathbf{S}_i[t]\mathbf{A}_{ii}$, which indicates that the next state for node $i$ only depends on its internal MC through $\mathbf{A}_{ii}$.
\end{proof}
\emph{Proposition III:} If for all the nodes $i$ and $j$ ($i,j\leq n$) we have $c_{ij}=1$ then the constraint-based model collapses to the original IM.
\begin{proof}
Similar to proof of the previous proposition, if for any node $i$ we set $c_{ij}=1$ for all $j\leq n$ in (\ref{EQ:expand}) then we get $\mathbf{p}_i[t+1]=\sum^n_{j=1}d_{ij}\mathbf{S}_i[t]\mathbf{A}_{ij}$, which is the statement based on the original IM for the probability of next step states for node $i$.
\end{proof}

In the constraint-based model, if we set $c_{ij} = 0$ then the influence that node $i$ receives from node $j$ gets deactivated and instead the influence of node $i$ from itself increases proportionally to compensate for the absent influence to preserve the row stochastic property of the influence matrix.

Addressing the second problem mentioned in Section~V is straight forward. The key idea is to allow time varying MCs for the nodes. In other words, unlike the original IM that assumes fixed transition probability matrix for the MC of each node (i.e., fixed $\mathbf{A}_{ii}$s), the transition probability matrices of internal MCs in the DCIM are functions of influences and may dynamically change over time. As such, a set of potential internal MCs is associated with each node while only one of them will be used at each time. The selection of which MC to use at each time depends on the influences and the state of the network.
To give an example of such selection function, consider the following scenario in the network of computing nodes. Depending on how many nodes a node can influence (can send workload to) a different MC is required to capture the internal state transitions of the node. For each node with the out-degree (defined as the number of out-going links in the influence network) equal to $k$ we define $k+1$ MCs. In general, if it is influencing $h$ neighbors ($0\leq h \leq k$), its MC will be represented by $\mathbf{A}_{ii}^{(h)}$. To define a function that selects the right MC for the node depending on the number of nodes it is influencing one can consider the following function
\begin{equation}
  \mathbf{A}_{ii}{(x)}=\sum_{h=0}^{k}\delta_x(h)\mathbf{A}_{ii}^{(h)},
\end{equation}
where $x$ represents the number of nodes being influenced by node $i$, $\delta_x (h)=1$ if $x=h$ and 0 otherwise. Notice that $x=\sum_{j=1}^n c_{ij}$.
%We assume that when a node influences $k$ nodes, it sends the same amount of workload to each.

\subsection{Steady-State Analysis of DCIM}
We briefly discussed the steady-state analysis of the original IM in Section II. It was explained that in IM $\lim_{t \to \infty} \mathbf{S}[0]\mathbf{H}^t$ converges to the state occupancy probability vector of the network if $\mathbf{H}$ has an eigenvalue at 1 that dominates all of its other eigenvalues and the steady-state vector is independent of the initial state.
In this section, we discuss the convergence properties of the DCIM to a steady state vector of probabilities. To investigate the convergence properties of $\lim_{t \to \infty} \mathbf{S}[0]\mathbf{H}^t$ for DCIM, we should note that matrix $\mathbf{H}$ in DCIM is no longer a fixed matrix in time and evolves dynamically over different steps. Dynamic nature of $\mathbf{H}$ has two sources as described next. (1) The conditional influence matrix, $\mathbf{C}$ is a function of the current state of the network (i.e., a function of $\mathbf{S}[t]$). Recall that matrix $\mathbf{E}$ depends on matrix $\mathbf{C}$ based on (\ref{EQ:E}) and $\mathbf{H}$ depends on $\mathbf{E}$ according to (\ref{EQ:NewH}) and thus matrix $\mathbf{H}$ is dynamic and changes as the state of the network changes. (2) The second source of dynamic nature of $\mathbf{H}$ is the dynamic changes of the internal MCs. At each time, the MC of node $i$ might be different from that of the previous state and this will change the $A_{ii}$ matrix, which appears in (\ref{EQ:NewH}).

In light of the above, the steady state analysis of the DCIM is more complicated due to the dynamic nature of $\mathbf{H}$ matrix.
Suppose $\sigma$ is the set of all $\mathbf{H}$ matrices for a given network modeled by DCIM. Notice that $\sigma$ is a finite set, because the number of $\mathbf{H}$ matrices in $\sigma$ is only a function of the number of constraint matrices and the number of internal state transition matrices for the nodes. Hence, it is finite and its members can be indexed by the set $I = \{i\in \mathbb{N} | i \leq k\}$, where $k$ is the cardinality of $\sigma$ (i.e., we have
$\sigma=\{\mathbf{H}_i| i\in \mathbf{I}\}$). An infinite sequence of $\mathbf{H}$ matrices are denoted by $K_r$ in which each matric is indexed by $h_j$, where $h_j\in \mathbf{I}$ and $j\in \mathbb{N}$. We also suppose that the set of all such sequences is denoted by $\Omega=\cup_r K_r$.
Note that $E(\mathbf{S}[t+1]|\mathbf{S}[0])$ depends on the chosen sequence $K_r$ of $\mathbf{H}$ matrices. In other words, for a given sequence $K_r$ in $\Omega$ until time $t$ we have
\begin{equation}
  E(\mathbf{S}[t+1]|\mathbf{S}[0])=\mathbf{S}[0]\mathbf{H}_{h_1}\mathbf{H}_{h_2}...\mathbf{H}_{h_t}.
\end{equation}
The existence of the steady-state probability vector for a network with given $K_r$ depends on whether or not the following limit exists:
\begin{equation}
  \mathbf{H}(K_r)=\lim_{t \to \infty}\mathbf{H}_{h_1}\mathbf{H}_{h_2}...\mathbf{H}_{h_t}.
  \label{EQ:Lim}
\end{equation}

The problem of convergence for multiplication of a set of matrices has been extensively studied in the field of linear algebra during the past two decades [13, 14]. It is shown in [13] that the limit of (\ref{EQ:Lim}) exists if the members of $\sigma$ have the following properties. (1) The matrices need to be right-convergent product set (RCP set). A set of matrices, i.e., $\mathbf{M}=\{\mathbf{A}_i|i\in \mathbb{N}\}$, is an RCP set if and only if for all $\mathbf{A}_i$s limit $a= \lim_{t \to \infty} (\mathbf{A}_i)^t$ exists.
(2) The matrices need to have the same left eigenspace for the eigenvalue 1 and their joint spectral radius is strictly smaller than 1.

As mentioned in Section~II.A and shown in [1], $\lim_{t \to \infty} (\mathbf{H}_i)^t$ exists only in certain conditions (i.e., when $\mathbf{H}$ has an eigenvalue 1 that dominates all of the other
eigenvalues). Similarly, in DCIM, $\lim_{t \to \infty} (\mathbf{H}_i)^t$ exists for all $i\in \mathbf{I}$ if they satisfy the mentioned condition. In that case, $\sigma$ would be an RCP set.
However, even in that case not every $\mathbf{H}_i$ matrix satisfies the second condition.
If the set $\sigma$  satisfies the two mentioned properties then the system will converge. However, unlike the steady-states of MCs and IM (that there is only one possible steady-state for the system), here, depending on the stochastic path of the system and for each random sequence $K_r$ in $\Omega$ we may get a different vector.

\section{Application of the DCIM to Stochastic Load Balancing in a Network of Computing Nodes}
Here, we extend the model for workload exchange in a network of computing nodes, presented in Section~II.~B, to a model, which captures the role of load-balancing policy rules using DCIM. We show that using the DCIM, one can investigate the performance of load balancing policies, where each policy specifies the rules of interaction and workload exchange among the nodes of a network of computing elements. Similarly, the DCIM can be applied to other complex network applications in which the nodes of the network interact based on certain rules such as in social networks.
\subsection{DCIM for Modeling Network of Computing Nodes}
The basics of the stochastic model for the network of computing nodes is presented in Section~II.~B. In this section, we consider a network of 30 computing nodes, which are connected over a random underlying topology. The effects of the underlying network topology have been studied in Subsection~V.~D. We assume that all the nodes have the same state space for their internal MC (i.e., underloaded, normal, and overloaded); however, the transition matrices for these MCs are different.
In order to capture the rule-based interactions based on the DCIM, we use the constraint matrix $\mathbf{C}$, as introduced in Section IV, which allows us to model different workload-balancing policies. To demonstrate how DCIM captures the effects of load-balancing policies, we first consider five plausible workload distribution policies:
\begin{itemize}
\item\emph{Policy 1:} Node $i$ gets influenced by (receives workload from) node $j$ if and only if (iff) node $i$ is underloaded and node $j$ is overloaded.
\item\emph{Policy 2:} Node $i$ gets influenced by node $j$ iff node $i$ is underloaded and node $j$ is overloaded or node $i$ is underloaded and node $j$ is in normal state.
\item\emph{Policy 3:} Node $i$ gets influenced by node $j$ iff node $i$ is underloaded and node $j$ is overloaded or node $i$ is in normal state and node $j$ is overloaded.
\item\emph{Policy 4:} Node $i$ gets influenced by node $j$ iff either node $i$ is underloaded and node $j$ is overloaded, node $i$ is underloaded and node $j$ is in normal state or node $i$ is in normal state and node $j$ is overloaded.
\item\emph{Policy 5:} Node $i$ gets influenced by node $j$ iff either node $i$ is underloaded and node $j$ is overloaded, node $i$ is underloaded and node $j$ is in normal state, node $i$ is in normal state and node $j$ is overloaded or node $i$ is in normal state and node $j$ is in normal state too.
\end{itemize}

For example, the constraint matrix $\mathbf{C}_3$ for Policy 3 is defined by $c_{ij}^3=s_{i3}s_{j1} + s_{i2}s_{j1}$, where $i \neq j$ and $c_{ij}^3$ denotes the entry in the $i$th row and $j$th column of $\mathbf{C}_3$ and
$s_{ab}$ denotes the binary variable corresponding to $b$-th (overload, normal, and undeload) possible state of node $a$.
The rest of the constraint matrices for other policies can be defined in a similar way using boolean logic.
Since influencing in this application means sending workload to other nodes, the influenced node can only stay in the normal state or go into overload state due to the influence. As such we consider the state-transitions for interactions among nodes, i.e., $\mathbf{A}_{ij}$, with a 3 by 3 matrix with rows equal to $[0.5, 0.5,  0]$.

\begin{figure*}
   	 \centering
    		\subfigure[]{
        		\includegraphics[width=2.9 in]{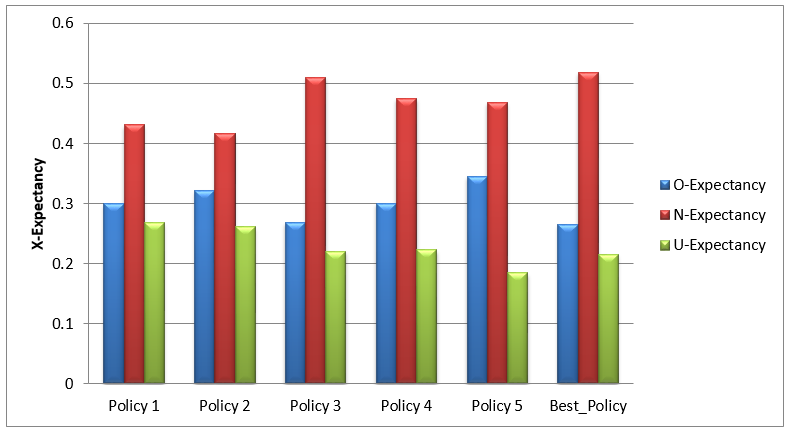}}
    	     \subfigure[]{
        		\includegraphics[width=2.5 in]{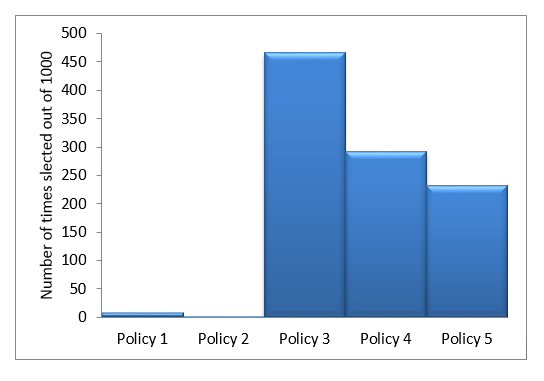} }
    	\caption{(a) Comparison of different load-distribution policies for a network with 30 computing nodes over 1000 steps and (b) number of times that the Best Policy approach chose each of the five policies over 1000 steps.}
    \label{FIG:Res1}
	\end{figure*}

%\begin{equation}
%\mathbf{A}_{ij}=
% \begin{pmatrix}
%  0.5 & 0.5 & 0  \\
%  0.5 & 0.5 & 0 \\
%  0.5 & 0.5 & 0
% \end{pmatrix}
%\end{equation}

\subsection{Performance Evaluation of Load-Balancing Policies}
One approach to quantify the performance of a load distribution policy is to calculate the expected normal-state occupancy ratio of the nodes if the goal is a balanced system. For other load-distribution goals, one can define various performance metrics. In general and independent of the selected performance measures, analyzing the performance of the policy using the steady-state analysis and asymptotic analysis of the stochastic model if possible would be the most effective approach. However, based on our discussion in Section~IV.~B about the steady-state behavior of DCIM, the convergence to a unique steady-state is not guaranteed to allow for such analysis. As such, in this section, we analyze the performance of the policies implemented based on DCIM using simulations.

To measure the performance of different load-distribution policies, we first apply the above policies to the network and for each policy we collect information about the internal state of the nodes during a sequence of 1000 steps and averaged over 1000 different initial states for the nodes.

%\begin{figure}[tbh!]
% 	\centering
% 	\includegraphics[width=3.5 in]{Result1-2}
% 	\caption{Comparison of different load-distribution policies for a network with 30 computing nodes over 1000 steps.}
% 	\label{FIG:Res1}
% \end{figure}

Given the current state of the network (internal state of all the nodes) and given the workload distribution policy, DCIM enables us to calculate the probability of occupying different states in the next time step in the vector $\mathbf{p}[t+1]$. Based on this vector, we define the notion of  {\it overall $X$-expectancy} as the expected portion of nodes being in state $X$ asymptotically over time. More precisely, $X$-expectancy is defined as:
\begin{equation}
  X_{\text{expc}}=\lim_{t \to \infty} \frac{\sum_t\sum_{i=1}^{n}\mathbf{p}_{i,X}[t]}{nt},
  \label{EQ:EXPC}
\end{equation}
where $\mathbf{p}_{i,X}[t]$ is the probability of node $i$ being in state $X$ at time $t$.
%\begin{figure}[tbh!]
% 	\centering
% 	\includegraphics[width=2.0 in]{Result1-2-1}
% 	\caption{Number of times that the Best Policy approach chose each of the five policies over 1000 steps.}
% 	\label{FIG:Res2}
% \end{figure}
Similar to this definition, we call the probability of being in a normal state for a computing node in the next time step, the {\it step-wise $N$-expectancy}. If a balanced computing network is desirable, one can choose a sequence of load-distribution policies that increases the overall $N$-expectancy. In case that increasing the overall $N$-expectancy is not feasible due to restrictions on asymptotic analysis of the model, one can find a sub-optimal sequence of policies by increasing step-wise $N$-expectancy. In Fig.~3-a, we show the overall $X$-expectancy values for overload, underload and normal states after 1000 steps in time for the five policies introduced earlier. Based on this result, we observe that load-distribution based on Policy 3 can lead to higher overall $N$-expectancy value (i.e., more normal states) and thus more balanced load distribution under similar conditions compared to other policies. We can also observe that some policies lead to higher $O$-expectancy values than others which means more overloaded nodes.

Given a list of plausible workload distribution policies and the current state of the network, $S[t]$, we can check to see which policy provides the highest step-wise $N$-expectancy if it was applied to the network at the next time step. We call the resulting policy, the {\it Best Policy} at time $t+1$. This means that a network can change load distribution policy in each step depending the status of the system to achieve better performance. The results shown in Fig.~\ref{FIG:Res1}-a compares the $X$-expectancy of the network after applying each of the policies mentioned earlier (fixed over time) in comparison to the Best Policy, which selects one of these five policies at each step depending on the state of the system. Figure~\ref{FIG:Res1}-b shows the number of times that each of the five policies was selected by the Best Policy approach over 1000 steps. In this example, we observe that the difference between the $N$-expectancy of the Best Policy and Policy 3 is not significant.

Next, we will use the above definitions to design the load distribution policy for the network, which leads to desirable status (in this case highest step-wise $N$-expectancy) at each step, without considering a given list of plausible policies.
%\begin{figure}[tbh!]
% 	\centering
% 	\includegraphics[width=3.0 in]{Result2}
% 	\caption{{\bf To be Updated}}
% 	\label{FIG:Res2}
% \end{figure}
\vspace{-0.3cm}
\subsection{Designing Optimum Load Distribution Policy}
Based on the discussion in the previous subsection one can use the $N$-expectancy to identify the policy among a given set of policies that leads to highest overall (in certain cases) or step-wise $N$-expectancy value. In this subsection, we address an even more interesting problem, which is how to design the optimum policy for our load distribution goal. In other words, we discuss designing the right constraint matrix $\mathbf{C}$ to achieve a specific load distribution goal for the next step.
\subsubsection{Algorithm for Optimum Policy Identification}
Identifying the right constraint matrix for a network of interacting nodes based on the global state of the system to achieve ceratin goal implies that instead of applying the same rules to all of the nodes, different nodes may interact differently depending on the effect of their local interactions on the global behavior.
For instance, if the goal of the load-distribution policy is balancing the load in the network, we need to identify matrix $\mathbf{C}$ that results in the largest $N$-expectancy in each step based on the current state of the system. The brute-force approach to identify such matrix is to search the space of all possible $\mathbf{C}$ matrices (i.e. check all possible workload distribution mechanisms) to find the one with the largest step-wise $N$-expectancy (Algorithm 1).
Next, we will discuss the complexity of this algorithm and identify conditions that allows one to find the optimum constraint matrix with less computational complexity than that of the brute-force approach presented in Algorithm 1.

\subsubsection{Complexity Analysis}
 Assuming that the network has $n$ nodes and the average degree of the nodes in the network is $k$ then the cardinality of $\mathcal{C}$ is in the order of $2^{kn}$. Therefore, the complexity of the Optimum Policy Selection algorithm presented in Algorithm 1 is simply $O(2^{n})$. Unfortunately, the exponential growth of the complexity of the algorithm makes it inapplicable to even moderate sized networks.
To understand under what conditions the complexity of this algorithm will be lower than exponential, let us look at (\ref{EQ:expand}) and expend it further for a specific state, say probability of node $i$ going to state normal in the next state as shown below:
 \begin{equation}
\begin{aligned}
\mathbf{p}^{\text{state }N}_i[t+1] &= (d_{ii}+\sum_{j \neq i}^n d_{ij}(1-c_{ij}))\mathbf{S}_i[t]\mathbf{A}_{ii}(x)\begin{pmatrix} 0 \\ 1\\0 \end{pmatrix}\\
&+ (\sum_{j\neq i}^n d_{ij}c_{ij})S_j[t]\mathbf{A}_{ji}\begin{pmatrix} 0 \\ 1\\0 \end{pmatrix}.
\end{aligned}
\label{EQ:expand2}
\end{equation}
Since $c_{ij}$s depend on $\mathbf{A}_{ij}(x)$s and vice versa, determining  $c_{ij}$s that maximizes $\mathbf{p}^{\text{state }N}_i[t+1]$ requires the complete search of $\mathcal{C}$.
However, under the assumption that the internal MCs of the nodes are fixed, we could find the constraint matrix $\mathbf{C}$ that provides the highest $N$-expectancy.
In the latter case, $\mathbf{p}^{\text{state }N}_i[t+1]$ does not depend on the constraint value of other nodes. In other words, the value of $\mathbf{p}^{\text{state }N}_i[t+1]$ only depends on the fixed parameters $\mathbf{D}$, $\mathbf{A}_{ij}$ and $\mathbf{S}[t]$. Furthermore, we have
\begin{equation}
\begin{aligned}
\mathbf{p}^{\text{state }N}[t+1]=\sum_{i=1}^n \mathbf{p}^{\text{state }N}_i[t+1].
\end{aligned}
\label{EQ:expand3}
\end{equation}
and hence for all $1\leq i\leq n$, $\mathbf{p}^{\text{state }N}_i[t+1]$ is positive then maximizing $\mathbf{p}^{\text{state }N}[t+1]$ for all the nodes in the network is equivalent with maximizing the individual $\mathbf{p}^{\text{state }N}_i[t+1]$. To do so, we simply need to compare the first and second terms in (\ref{EQ:expand2}) and set $c_{ij}=0$ if the first term is larger and set $c_{ij}=1$ otherwise. These steps are presented in Algorithm 2, with $O(kn)$ complexity.

\vspace{-0.3cm}
\begin{algorithm}[H]
\caption{Optimum Policy Selection Algorithm}
\label{CHalgorithm}
\hspace*{-0.7cm}
\begin{algorithmic}[1]
\State $\mathcal{G}$ \Comment A directed underlying topology for the network of interactions (e.g., topology of the communication network for transferring workloads).
\State $\mathbf{A}_{ij}(x)$ \Comment State-transition matrices of network nodes.
\State Identify set $\mathcal{C}$, which is the set of all possible constraint matrices $\mathbf{C}$ based on $\mathcal{G}$.
\State Initialize $C_{\text{opt}}$ with a random $\mathbf{C} \in \mathcal{C}$.
\State Calculate $N^{\text{opt}}_{\text{expc}}$ for $C_{\text{opt}}$ using $\mathbf{p}[t+1]$ \Comment Calculating the step-wise $N$-expectancy for $C_{\text{opt}}$ Policy.
\For    {All $\mathbf{C} \in \mathcal{C}$}
 \State Calculate the $N_{\text{expc}}$ for the current $\mathbf{C}$ using $\mathbf{p}[t+1]$.
\If {$N_{\text{expc}}>N^{\text{opt}}_{\text{expc}}$}
\State Replace the old $C_{\text{opt}}$ and $N^{\text{opt}}_{\text{expc}}$ by the current $\mathbf{C}$ and $N_{\text{expc}}$, respectively.
\EndIf
\EndFor
\State Return $C_{\text{opt}}$. \Comment $C_{\text{opt}}$ is the identified optimum load distribution rule for the network based on its current state.
\end{algorithmic}
\end{algorithm}
\vspace{-0.3cm}

\begin{figure*}
   	 \centering
    		\subfigure[]{
        		\includegraphics[width=2.8 in]{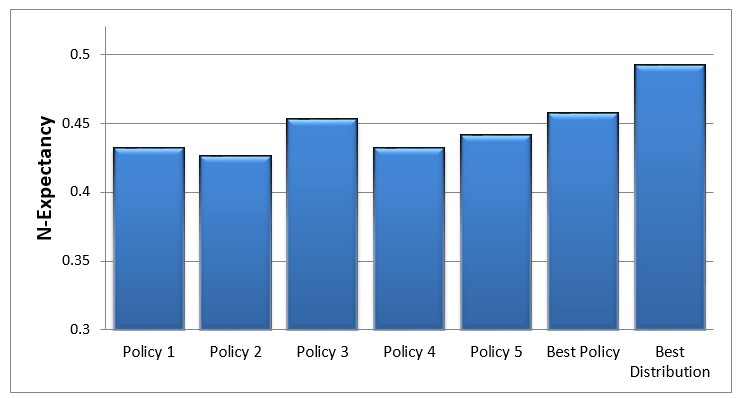}}
    	       \subfigure[]{
        		\includegraphics[width=2.8 in]{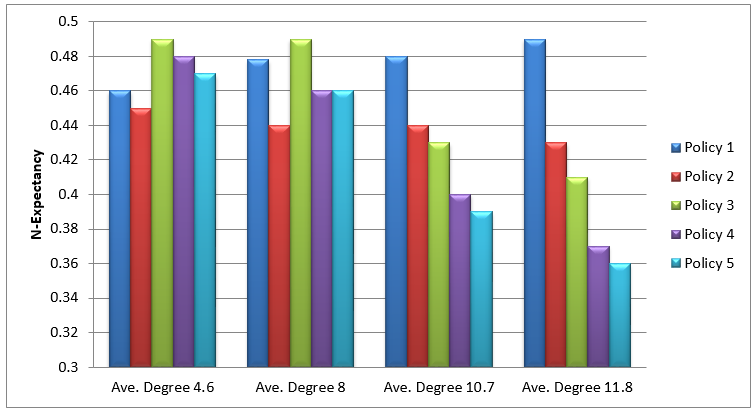}}
    	\caption{(a) Performance of the optimum load distribution policy in comparison with the rest of the mechanisms, and (b) the $N$-expectancy of the computing network for the same policies but different topologies.}
    \label{FIG:Res3}
	\end{figure*}

In Fig.~\ref{FIG:Res3}-a, we have presented the overall $N$-expectancy for the optimum constraint matrix along with the five policies introduced earlier and the Best Policy approach. We can observe that the designed constraint matrix or policy can lead to better load balancing performance compared to the rest of the load distribution mechanisms.

% \vspace{-0.1in}
%\begin{figure}[tbh!]
% 	\centering
% 	\includegraphics[width=3.0 in]{Result1-2-2}
% 	\caption{Performance of the optimum load distribution policy in comparison with the rest of the mechanisms.}
% 	\label{FIG:Res3}
% \end{figure}
%  \vspace{-0.2in}

\subsection{Impacts of the Underlying Topologies}
In the previous subsections, we assumed that the topology of the network is given and fixed. In this section, we show that the underlying topology of the network also affects the interaction behavior. We have specifically studied the impact of various topologies with various degree of connectivity when the load distribution policies introduced earlier are applied to the networks, on the network load distribution performance. Figure~\ref{FIG:Res3}-a represents the results for overall $N$-expectancy for various networks with different degree of connectivity. In this example, degree specifies the logical connection among the computing nodes for exchange of the workload. The presented results imply that the connectivity and influence relations among the nodes of the network impact the performance of the load distribution policy. For instance, Fig.~\ref{FIG:Res3}-b suggests that as the degree increases, Policy 1 does better than the other policies. The study of how exactly the topologies with different properties affect the interactions and behavior of the system is a future direction of this research. Identifying the optimum underlying topology for a network is also a very important problem to be investigated.

\vspace{-0.4cm}
\begin{algorithm}[H]
\caption{Optimum Policy Selection Algorithm with Fixed Internal MCs}
\label{CHalgorithm}
\hspace*{-0.7cm}
\begin{algorithmic}[1]
\State $\mathcal{G}$ \Comment A directed underlying topology for the network of interactions (e.g., topology of the communication network for transferring workloads).
\State $\mathbf{A}_{ij}$ \Comment Fixed state-transition matrices for network nodes.
\State $C_{\text{opt}}=\mathbf{0}$.
\For   {$i=1$ to $n$}
\For {nodes $j$ that are neighbor of node $i$ based on $\mathcal{G}$}
\If {$\mathbf{S}_i[t] \mathbf{A}_{ii}(x)\begin{pmatrix} 0 \\ 1\\0 \end{pmatrix} \leq \mathbf{S}_j[t] \mathbf{A}_{ji}\begin{pmatrix} 0 \\ 1\\0 \end{pmatrix}$}
 \State $c^{\text{opt}}_{ij}=1$.
 \Else
 \State $c^{\text{opt}}_{ij}=0$.
\EndIf
\EndFor
\EndFor
\State Return $C_{\text{opt}}$. \Comment $C_{\text{opt}}$ is the identified optimum policy for the network.
\end{algorithmic}
\end{algorithm}
\vspace{-0.4cm}

%\begin{figure}[tbh!]
% 	\centering
% 	\includegraphics[width=3.5 in]{Result1-2-3}
% 	\caption{The $N$-expectancy of the computing network for the same policies but different topologies.}
% 	\label{FIG:Res3}
% \end{figure}
% \vspace{-0.2in}

\section{Related work}
In this section, we briefly review the related work in two categories: (1) interaction and propagation models in networks,  and (2) applications of IM. We will also discuss the position of the current work relative to the works in these categories.

\emph{a. Interaction and Propagation Models in Networks:}
Network interactions and propagation models are pivotal problems in the network theory and complex systems and have been extensively studied in the literature. Examples of network interaction models besides the Influence model include the voter model, the contact process, the Ising model [15, 16], and the stochastic automata network model [17, 18]. Voter model, the contact model and the Ising model, which are usually referred to as interacting particle systems, are generally confined to infinite lattice structures. The basic idea of these models is that each node also has a clock rate and whenever the clock strikes, a coin flips and determines the next state of the node. The clock rate of a node depends on the number of its neighbors being in different states.
The stochastic automata network model, which perhaps is the closest theory to IM, considers networks consisting of interactive components modeled as automata [17, 18]. The model consists of all of the individual automata and the interactions of them over the network. In this model, a master MC is constructed to model the evolution of the whole system via a large transition matrix called the generator matrix. However, the major drawback of this approach is that the generator matrix's size increases exponentially as the number of the automata increases, and therefore, the theory suffers from complexity issues. Addressing this issue was one of the main motivations of the IM, which provides a computationally tractable framework for interactions.
Besides the models mentioned above, which focus on the detailed interactions of network components, analysis of propagation phenomena in the networks have been studies by percolation theory [19], sand-pile model (aka Bak-Tang-Wiesenfeld model) [20], SIR, SIS and SEIS models [21]. These models have been designed to capture simple interactions of the propagation process in usually large scale networks. In latter models, the details of dynamics of individual nodes in the network are not captured and the interactions are mainly based on the topology. In addition to these theoretical interaction models, problem specific data-driven interaction models have also been proposed. For instance, the interaction models for the problem of cascading failures in power systems have been reviewed in [4].

\emph{b. Applications of IM:}
IM has been applied to a variety of interesting real world problems. In [5] and [6], IM was used for modeling cascade of failures in power-grid
networks using the binary and also the evil rain versions of the IM.
In [22], IM is used for modeling an interdependent network to identify optimum inter-network service assignments in order to have a cascade-resilient network.  In a series of recent works [23-26], IM has been used for analyzing the interactions among the components during cascading failures in power grids. These works have a data-driven approach in using IM and use IM to infer the interactions among the components of the system during cascades in power grids. Another example of the application of IM is in modeling disease propagation [27]. In the latter paper, a combination of IM and the SIS model are used to model disease epidemics. Finally, several projects based on IM, such as modeling human interaction [28], automatic speaker identification (for documenting meetings audio recording) [29], and multi-sensor data fusion [30] as well as modeling opinion dynamics [31, 32] have been developed.
The application of the IM to distributed computing environments based on the rules of interactions is one of the novel contributions of the current paper.
\section{Conclusions}
Modeling interactions and influences between components of a system has been vastly studied in the literature and many models have been proposed; however, there are still many real world scenarios that do not quite fit into any of the available models. Examples of such systems are networks of computing nodes, which interact by exchanging workload based on load distribution policies. In this paper, we discussed that the original Influence Model has limitations in modeling dynamic and constraint-based interactions in the network. As such, we presented an extension of the IM, namely Dynamic and Constraint-based Influence Model (DCIM), to alleviate these limitations. We discussed the steady-state analysis for the DCIM and mentioned that the convergence properties of IM is not always valid for DCIM; however, one can still benefit from the application of the DCIM in modeling and analyzing the behavior of a network of interacting elements. We introduced the novel application of the DCIM to model a network of computing nodes for load balancing purposes and discussed designing optimum load-balancing policy based on DCIM. The proposed model in this paper enables the application of the Influence-based interaction modeling to more general network interaction scenarios.

\end{document}